\documentclass{article}
\usepackage{amssymb, amsmath}
\usepackage[all]{xy}
\usepackage[a4paper, textwidth=160mm, marginratio=1:1]{geometry}

\newtheorem{theor}{Theorem}[section]

\newtheorem{cor}[theor]{Corollary}

\newtheorem{lem}[theor]{Lemma}

\newtheorem{proposition}[theor]{Proposition}

\newtheorem{conjecture}[theor]{Conjecture}

\newtheorem{defin}[theor]{Definition}

\newtheorem{rem}[theor]{Remark}

\def\qed\hfill{$\squares$}

\def\build#1_#2^#3{\mathrel{\mathop{\kern0pt#1}\limits_{#2}^{#3}}}

\def\sm{\smallskip}
\def\m {\sm}

\def\lra{\longrightarrow}

\def\zmoddeux{\mathbb{Z}/2\mathbb{Z}}

\def\fdeux{\mathbb{F}_2}

\def\zdeux{{\mathbb{Z}/2\mathbb{Z}}}

\def\fd{f\mathcal{D}}

\def\fdn{f\mathcal{D}_n}

\def\fd{f\mathcal{D}}

\def\fdn{f\mathcal{D}_n}
\def\fdnp1{\fd _{n+1}}

\def\etendnx{\mathbf{E}_n (X)}
\def\etend{\mathbf{E}}

\def\etendnnz{\mathbf{E}_{n-1} (Z)}
\def\etennd{\mathbf{E}_{n-1}}
\def\etenrd{\mathbf{E}_{r}}
\def\rp{\mathbb{R}P}

\def\sodeux{SO (2)}

\def\shuffle{\mathrm{sh}}
\newcommand{\J}{\mathbb{J}}

\def\hur{\mathcal{H}}

\def\act{\mathfrak{a}}

\begin{document}
\title{Relations in the modulo $2$ homology of framed disks algebras}
\author{G\'erald Gaudens }
\maketitle

\begin{abstract}
We study in this work the homology structure of spaces having
an action of the \emph{framed disks} operad. In particular, we
compute the relations between Kudo-Araki operations and
generalized Batalin-Vilkovisky operators. As an application, we
complete the computations of the modulo $2$ homology of $\Omega^2
S^3$ as a Batalin-Vilkovisky algebra, and give some evidence for a
classical conjecture about the modulo $2$ Hurewicz homomorphism of
the infinite loop space of the sphere spectrum.
\end{abstract}

\section{Introduction}
\label{intro}


\subsection{Summary of results}
\label{summary}

Let $X$ be a well pointed space having the homotopy type of a
CW-complex, and consider the iterated loop space $\Omega^n X$ for
$n\geq 2$, and let $\mathbb{K}$ be a ground coefficient field for
the singular homology functor $H_*$. The natural  action of the
little $n$-disks operad $\mathcal{D}_n$ on $\Omega^n X$ provides
$H_*\Omega^n X$ with the structure of a
$H_*\mathcal{D}_n$-algebra. F.
Cohen~\cite{Cohen-Lada-May:homiterloopspaces} has shown that this
gives $H_* \Omega^n X$ the structure of an $e_n$-algebra (see
Definition \ref{definition algebre de Gerstenhaber} in the
appendix). Over the field of rational numbers $\mathbb{Q}$, these
two notions essentially coincide.

\sm Suppose now that $X$ is equipped with a (continuous) pointed action of
$SO(n)$, which means that $SO(n)$ acts on $X$ through pointed maps.
Getzler~\cite{Getzler:BVAlg}, Salvatore and
Wahl~\cite{Salvatore-Wahl:FrameddoBVa} have noticed that the
little $n$-disks operad $\mathcal{D}_n$ is an $SO(n)$-operad which
acts in the category of $SO(n)$-spaces on $\Omega^n X$, or
equivalently that the framed $n$-disks operad $f\mathcal{D}_n$
acts on $\Omega^n X$. Therefore $H_*\Omega^n X$ is an
$H_*SO(n)$-algebra over the $H_*SO(n)$-operad $H_*\mathcal{D}_n$, or in other words
$H_*\Omega^n X$ is an algebra over the operad $H_*f\mathcal{D}_n$.
One goal of this paper is to provide some explicit computations of
this structure, particularly the case  of $H_*\Omega^2 S^3$ left
unachieved in \cite{gaudensmeni0}.

\sm Obviously, the structure of $H_*f\mathcal{D}_n$-algebra is the
structure of an $e_n$-algebra together with the structure of
$H_*SO(n)$-module which satisfy some compatibility relations
between them. Restricting to the field of rational numbers $\mathbb{Q}$ for $n=2$, Getzler~\cite{Getzler:BVAlg} has shown that
$H_*(f\mathcal{D}_2;\mathbb{Q})$-algebras are exactly
$BV_2$-algebras (i.e. Batalin-Vilkovisky algebras, see Appendix \ref{previous}). Over fields of
positive characteristic, $H_*(f\mathcal{D}_2;\mathbb{K})$-algebras
are still $BV_2$-algebras, but carry Kudo-Araki operations as
well. More generally, the structure of $H_*f\mathcal{D}_n$-algebras over the field
$\mathbb{Q}$ was explored in \cite{Salvatore-Wahl:FrameddoBVa,gaudensmeni0}. We will study in this work this structure for $\mathbb{K}= \mathbb{F}_2$.

\sm As an example of a typical relation (Corollary
\ref{comm:bv:q}) we get that the Batalin-Vilkovisky operator
$\Delta$ and the Kudo-Araki operation $Q_1$ in the homology of a
two fold loop space satisfy for all $x$ in $H_*\Omega^2 X$
$$
\Delta (Q_1 x)= \{\Delta (x), x\}+ \Delta (x)* \Delta
(x)= \Delta (x*(\Delta x))
$$
where $*$ denotes the Pontryagin product.

\sm As explained above, this relation yields a complete computation of $H^*\Omega^2 S^3$ as a $BV$-algebra.

\sm {\bf Theorem \ref{bvcalc:omegasdeux}}.~~ {\it The action of
the Batalin-Vilkovisky operator $BV$ on $H_* (\Omega^2S^3 ,
\mathbb{F}_2 )= \mathbb{F}_2 [u_i; i\geq 1]$ is  given on
monomials $u_{i_1}^{\ell_1}*u_{i_2}^{\ell_2}*\ldots
*u_{i_n}^{\ell_n}$ with $\ell_j>0$ by
\begin{eqnarray}
{BV} (u_{i_1}^{\ell_1}*u_{i_2}^{\ell_2}*\ldots *u_{i_n}^{\ell_n})= \sum_{j=1}^{n} \ell_j
u_1^{2^{i_j}}* u_{i_1}^{\ell_1}*u_{i_2}^{\ell_2}*\ldots *u_{i_j}^{\ell_j-1}*  \ldots * u_{i_n}^{\ell_n}.
\end{eqnarray}
}

\sm Recall from \cite{Salvatore-Wahl:FrameddoBVa} that an
$(n-1)$-connected group-like $f\mathcal{D}_n$-algebra is  weakly
equivalent as an $f\mathcal{D}_n$-algebra  to an $n$-fold loop
space on a certain pointed $SO(n)$-space. It is therefore a mild
restriction to think of $f\mathcal{D}_n$-algebras as $n$-fold loop
spaces on pointed $SO(n)$-spaces, although our results hold for
general $f\mathcal{D}_n$-algebras.  In \cite[section
3]{gaudensmeni0}, we explained why for computational purposes, it
is further a small restriction to focus on
$f\mathcal{D}_n$-algebras that are $n$-fold loop spaces $\Omega^n
X$ over a space  $X$ with trivial $SO(n)$-action. In this case,
the global action of $SO(n)$ on $\Omega^n X$ is \emph{not}
trivial, and in fact the action of $H_* SO(n)$ on $H_*\Omega^n X$
is in this case nicely related to the $J$-homomorphism. This
relation was central to \cite[section 5]{gaudensmeni0}. This fact
and the computation  in \cite[section 6]{gaudensmeni0} (see also
\cite{menibv})  will serve as basis for the computations in
Section \ref{comp:hur:qso} and Section \ref{sec:bvspheres} of the
present work. We prove:

\sm {\bf Theorem \ref{divtrivial}}.~~ {\it All classes in
$\pi_{>0} QS^0$ divisible by the image of the $J$-homomorphism are
annihilated by the {\it mod} $2$ Hurewicz homomorphism, except the
Hopf maps and their composition squares (which are \emph{not}
annihilated by the Hurewicz homomorphism). In other words,
Conjecture \ref{curtismadsen} holds on the ideal generated by the
image of the $J$-homomorphism. }

\sm In presenting this result on the Hurewicz homomorphism for
$QS^0$, we do not claim originality. Although we are not aware of
any published reference, it is difficult to believe that such a
result is not known to any expert. Rather, these computations are
meant to advertise that group operations (the group being in our
case $SO(n)$) on spaces might help to analyze the Hurewicz homomorphism.
More work in this direction appears in \cite{gaudensmeni2}. We
also mention that it is possible to push a little bit further the
idea of analyzing the Hurewicz homomorphism for $QS^0$ by higher
composition methods ({\it i.e.} Toda brackets).

\subsection{Organisation of the paper}
\label{organisation}

From now on we assume that $\mathbb{K}= \mathbb{F}_2$. After the
necessary recollections on Kudo-Araki operations (Section
\ref{recoll:dyer:lashoff}), we explore the structure of
$H_*(f\mathcal{D}_n)$-algebras (Theorem \ref{theo:bat},
\ref{theo:bat:2}, and \ref{higherbv}). The
$H_*(f\mathcal{D}_n)$-structure allows one to produce new formulas
about composition products with elements coming from the
$J$-homomorphism in the homology of loop spaces.  We then compute
the action of $H_*(SO(2))$ on $H_*(\Omega^2 S^k)$ for all $k\geq
1$ (Theorem \ref{bvcalc:omegasdeux}, \ref{omegadeuxsdeux} and
\ref{trivialahah}), which was begun in \cite{gaudensmeni0}. In
Section \ref{comp:hur:qso}, we show how the methods of
\cite{gaudensmeni0} and of the present paper apply to give some
information about the Hurewicz homomorphism for $QS^0$, the
infinite loop space associated to the stable homotopy of spheres
(Theorem \ref{divtrivial}).

\subsection{Notations and conventions}

Recall that we assume all topological spaces to have the homotopy
type of a $CW$-complex. If $X$ is a topological space, the
singular homology groups are denoted by $H_*X$ and  are always
taken with $\zmoddeux$-coefficients; $D_X: X\lra X\times X$ is the
diagonal map of $X$; we usually, alleviate the notation by
omitting $X$ from it. For any continuous map $\Psi$, we let
$\Psi_*$
be its induced map in homology.
The map $\hur : \pi_* X \to H_* X$ is the {\it modulo} $2$
Hurewicz homomorphism.

\sm We fix  a non zero positive integer $n>0$. The bottom and top
homology classes in $H_* S^n$ are denoted respectively by $b_0$
and $b_n$. $S^n$ equipped with the antipodal action $T$ becomes a
$\zmoddeux$-space. We provide $S^n$ with its usual $\zmoddeux$
equivariant cell structure (with two antipodal cells  $e^i$ and
$Te^i$ in each dimension $0\leq i\leq n $).

\sm Given a cell complex $X$ (i.e a CW-complex together with a
specific cell decomposition), we let $C_*^{\mathrm{cell}} X$ be
the associated cellular chain complex. In  particular
$C_*^{\mathrm{cell}} S^n$ is the cellular chain complex of $S^n$
with respect to the above classical $\zmoddeux$-cell structure,
and is denoted by $W_\bullet$, and is the standard periodic
resolution of $\zmoddeux$ over $\fdeux$ truncated at level $n$:
$$
(W_n)_i= \fdeux e^i \oplus \fdeux Te^i~~, ~~ \mathrm{for}~ i\leq n
$$
and $d(e^i)= e^{i-1}+Te^{i-1}$, and $(W_i)=0$ for $i>n$.

\sm The operad of little $n$-disks consists of the  symmetric
sequence of spaces $(\mathcal{D}_n(i))_{i\geq 0}$ of TD-embeddings
({\it i.e.} affine embeddings that are compositions of
translations and dilatations) of the disjoint union of $i$-little
$n$-discs in the interior of a standard $n$-disc, whose images are
pairwise disjoint, equipped with with the suitable topology. The
symmetric group on $i$ letters acts on $\mathcal{D}_n(i)$ in the
obvious way, and $\{\mathcal{D}_n(i)\}_{i\geq 0}$ forms a
topological operad by ``\emph{putting big disks in small embedded
discs}''. Any $n$-fold loop spaces is an algebra over this operad,
and it happens that any connected algebra over the little
$n$-discs operad is actually equivalent as a
$\mathcal{D}_n$-algebras to an $n$-fold loop space. The reference
for this material is the classical \cite{MayJ:Geoils}.


\sm Recall that the space of little $2$-disks in $D^n$ is homotopy equivalent
to $S^{n-1}$, equivariantly with respect to the actions of
$\zmoddeux$ and $SO(n)$.

\sm The little $n$-discs operad is a sub-operad of the framed
little $n$-discs operad $\fdn$ which consists of the semi-direct
product of $\mathcal{D}_n$ with the linear group $SO(n)$ (see
\cite{Salvatore-Wahl:FrameddoBVa}).

\sm The symbol $\delta_i^j$ for $i, j$ integers is the Kronecker
symbol, whose value is $1$ if $i=j$ and $0$ otherwise.

\sm If $X$ is an H-space, we let $*$ denote the Pontryagin product
associated to the H-space structure on $H_* X$.

\sm {\it Acknowledgements.-} We wish to thank B. Fresse, S. Kallel
and V. Tourchine for useful discussions. This project is a follow
up of the joint work \cite{gaudensmeni2} with L. Menichi,  whose
remarks have been invaluable to the author. This work has been
elaborated in the topology group of the University of Bonn, where
the author has found great support and friendship; it was part of
the author's {\it Habilitationsschrift}. The author also thanks
the Max Planck Institut f\"ur Mathematik of Bonn for supporting
the author while this work has been finished.

\section{Recollection on Kudo-Araki operations}
\label{recoll:dyer:lashoff}

\subsection{Extended products and their homology}
\label{otherview}

\subsubsection{Extended products.}
\label{natur}

The classical Kudo-Araki operations form a set of {\it modulo} $2$
homology operations for $\mathcal{D}_n$-algebras
\cite{kudo2,kudo1}. We give here a short recollection.

\sm
Let $n>0$. $\zdeux$ acts on $S^n$ by the antipodal action,  and on
$X\times X$ by exchanging factors. Let $\etendnx = S^n\times
_\zdeux (X\times X)$. We have a fiber bundle with structural group
$\zdeux$
$$
X\times X \lra \etendnx \lra \rp ^n
$$
and a Serre spectral sequence
$$
H_*^{\mathrm{tw}} (\rp ^n, H_* X^{\otimes 2}) \Rightarrow H_*
(\etendnx)
$$
where $H_*^{\mathrm{tw}}$ stands for homology with twisted
coefficients. It is classical (see for instance \cite{Maygeneral,milgramstable,zaratithese}) that this spectral sequence collapses
at the $E_2$ term. Hence there is a vector space isomorphism
$$
H_{\bullet} \etendnx \cong H_\bullet^{\mathrm{tw}} (\rp ^n,
H_\bullet X ^{\otimes 2})~~.
$$
There are  classical isomorphisms \cite[lemma 1.1]{Maygeneral}:
\begin{eqnarray}
\label{isos} H_*^{\mathrm{tw}} ( \rp ^n, H_* (X\times X))&=&  H_*
(C^{\mathrm{cell}}_* S^n \otimes_{\zmoddeux} H_* (X^{\times 2}))
\\
\nonumber &=& H_* (W_n\otimes_{\zmoddeux} H_* (X^{\times 2}))
\\
\nonumber &\cong& H_* (W_n\otimes_{\zmoddeux} H_* (X)^{\otimes
2})
\end{eqnarray}
More generally, let $A$ be a \emph{free} $\zmoddeux$-CW-complex
and let $X$ be any space. There is a fibration sequence as before:
$$
X\times X \to A\times_{\zdeux} (X \times X) \to (A)_{\zdeux}
$$
We have isomorphisms:
\begin{eqnarray}
\label{isos2} H_*^{\mathrm{tw}} ((A)_{\zdeux}, H_* (X\times
X))&\cong& H_*
(C^{\mathrm{cell}}_* A \otimes_{\zmoddeux} H_* (X^{\times 2}))\\
\nonumber &\cong& H_* (C^{\mathrm{cell}}_* A\otimes_{\zmoddeux} H_* (X)^{\otimes
2})
\end{eqnarray}
These isomorphisms are natural in the variable $A$ with respect to equivariant
cellular maps and in $X$ with respect to all maps.

\sm Furthermore, if $(g_i)_{i\in I}$ is a basis (in particular the
indexing set $I$ is totally ordered) for $H_* X$, a basis for the
twisted homology $H_*^{\mathrm{tw}} ( \rp ^n, H_* X^{\otimes
2})\cong H_* \etendnx $ is given by the list of homology classes
\begin{eqnarray}
\label{basis} \left[ e^0 \otimes g_i\otimes g_j \right] ~~,~~ i<j;
~~\left[e^n \otimes g_i\otimes g_j\right]~~, ~~i<j ~~;~~ \left[e^k
\otimes g_i \otimes g_i\right] ~~,~~0\leq k \leq n \nonumber&
\end{eqnarray}
for $i, j \in I$.

\subsubsection{Operations on the homology of $\mathcal{D}_{n+1}$-algebras.}

Let $X$ be a $\mathcal{D}_{n+1}$-algebra; $X$ is most typically an
$(n+1)$-fold loop space. There is an operadic action map:
$$
\theta: S^{n}\times X\times X \lra X
$$
which factorizes as
$$
\bar\theta : \etendnx \lra X ~~.
$$
We define:
\begin{eqnarray}
\theta_* (b_0 \otimes x\otimes y)&=& \bar\theta _* ([e^0 \otimes x\otimes y])  = x*y\\
\bar\theta_* ([e^i \otimes x\otimes x]) &=& Q_i x \mathrm{~for~}0\leq i\leq n\nonumber\\
\theta_* (b_n \otimes x\otimes y)&=&\bar\theta ([e^n\otimes x\otimes  y]) = \{x, y\}\nonumber
\end{eqnarray}
\sm
We have used the following  classical computation of the map $\zeta_*$ induced by
the quotient map
$$\zeta: S^n \times X \times X \lra \etendnx$$
\begin{eqnarray}
\label{quotientmapcomput}
\zeta(b_0 \otimes g_i \otimes g_j)&=&  [e^0 \otimes g_i \otimes g_j]~~, \mathrm{~for~} i\neq j\\
\zeta(b_n \otimes g_i \otimes g_j)&=&  [e^n \otimes g_i \otimes
g_j]~~, \mathrm{~for~} i\neq j \nonumber
\end{eqnarray}
and $\zeta_*$ is zero otherwise.

\sm Elements of the form $[e_i \otimes x\otimes x]$ are called
Kudo-Araki elements, those of the form $[e^0 \otimes x\otimes
y]$ are called Pontryagin elements, and those of the form
$[e^n\otimes x\otimes y]$ for $x\neq y$ are called brackets. Of
course $*$ is the Pontryagin product, $\{-,-\}$ is the Browder
bracket, and the operators $Q_i$ are the Kudo-Araki operations.
One checks that this definition and others are equivalent.

\sm The natural quotient map $S^n \lra S^n \vee S^n$ obtained by
collapsing the $(n-1)$-skeleton is $\zmoddeux$-equivariant and
induces a natural quotient map
\begin{eqnarray*} \xi:\etendnx
&\lra& S^n \wedge  (X\times X)
\end{eqnarray*}
that will be crucial to us in order to elucidate the structure of
$H_* (f\mathcal{D}_{n+1})$-algebras. We will need to know that
this quotient map satisfies \cite[part I, chap 3, p.
37-39]{milgramstable}:
\begin{eqnarray}
\label{sacre:formules}
\xi_* ([e_0\otimes x\otimes y])&=&0 \\
\nonumber \xi_* ([e_i\otimes x\otimes x])&=& 0 \mathrm{~for ~} i\leq n-1\\
\nonumber \xi_* ([e_n\otimes x\otimes x])&=& b_n \otimes x\otimes x\\
\nonumber  \xi_* ([e_n\otimes x\otimes y])&=& b_n\otimes x\otimes
y+ b_n \otimes y \otimes x
\end{eqnarray}
Observe that in\cite{milgramstable}, the author works with cohomology; we leave it to the reader
to make the necessary changes to get this result.

\begin{rem}
\label{renorm} One usually defines `upper index' Kudo-Araki
operations by
$$
Q^i x := Q_{i-|x|} x\quad .
$$
In this way, the operation $Q^i$ simply raises the degree by $i$.
These reindexed Kudo-Araki operations satisfy \emph{unstability}:
$Q^i X=0$ if $i<|x|$.
\end{rem}

\subsection{Properties of the Browder bracket and the Kudo-Araki operations}

The following holds for all little $(n+1)$-disks algebras, hence
in particular for \emph{framed little $(n+1)$-disks algebras},
{\it i.e.} algebras over $\fdnp1$.

\sm\begin{proposition}\cite[Theorem 1.3 (5) p.
218]{Cohen-Lada-May:homiterloopspaces} \label{quadratic:formula}
On the homology of a $\mathcal{D}_{n+1}$-algebra $X$, the
operations $Q_i$ are additive for $i<n$ and  $Q_{n}$ is quadratic
with respect to the Browder bracket, which means that for all $x, y \in
H^* X$
\begin{eqnarray}
Q_{n}(x+y)= Q_{n}(x)+ Q_{n}(y) +\{x,y\} \quad .
\end{eqnarray}
\end{proposition}

\begin{cor}
The Kudo-Araki operations determine Browder brackets in the modulo
$2$ homology of $\mathcal{D}_{n+1}$-spaces.
\end{cor}

In order to elucidate the structure of $H_*
(f\mathcal{D}_{n+1})$-algebra at the prime $2$, it is therefore
crucial to understand the relation between the $H_*SO
(n+1)$-action with Pontryagin products and Kudo-Araki operations.
We do this in Section \ref{bv:dl}. A direct consequence of
Proposition \ref{quadratic:formula} is:
\begin{rem}
\label{rem:one:one}\cite[Theorem 1.2 (3) p.
215]{Cohen-Lada-May:homiterloopspaces}  In the {\it modulo} $2$
homology of a little (n+1)-disks algebra $X$, there is a relation
for all $x\in H_* X$:
\begin{eqnarray}
\label{crochetnuleq} \{x,x\}=0 .
\end{eqnarray}
\end{rem}
\begin{rem}
\label{rem:three:three} Another fact worth pointing out is the
fact that the Browder bracket on the homology of an $n$-fold loop
space vanishes as soon as the $n$-fold loop space structure
extends to an $(n+1)$-fold loop space structure. In this case, a
Browder bracket associated to the $(n+1)$-fold loop space
structure is defined and the latter does not need to vanish.
\end{rem}

\section{The structure of  the {\it modulo} $2$ homology of $H_* f\mathcal{D}_n$-algebras}
\label{dl:et:bv}

Let $n\geq2$ and let $Z$ be an $f\mathcal{D}_n$-algebra. In this
section, we study the relations between the $H_* SO (n)$ action on
$H_* Z$ due to the $SO(n)$-action coming from $f\mathcal{D}_n$ and
the mod $2$ Kudo-Araki operations due to the action of the operad
of little $n$-disks in the operad $f\mathcal{D}_n$. We first need to
recall some elementary facts about the homology of $SO(n)$ in
order to proceed to the computations.

\subsection{Recollection on the homology of orthogonal groups}

We consider the reflection map
$\kappa:\mathbb{R}P^{n-1}\longrightarrow SO(n)$, that sends a line
to the reflection about its orthogonal hyperplane, then composed
with a fixed reflection (in order to get an orientation preserving linear transformation). Recall that
$$
H_* SO(n)= \Lambda (d_1, d_2,\ldots, d_{n-1})$$ where $d_{i}:=
\kappa_* e_i$ is the image of the unique non-zero class  $e_i \in
H_i \mathbb{R}P^{n-1}$ by $\kappa_*$. In particular, the diagonal
applied on the generators is \emph{the Cartan diagonal}
$\mathrm{D}_*$:
$$
\mathrm{D}_* d_k =\sum_{i+j=k} d_i \otimes d_j \quad .
$$
A reference for this classical fact is \cite{semc}.
\begin{defin}
Let $Z$ be an $f\mathcal{D}_n$-space. In what follows, we denote
the action of $d_i\in H_*Z$ on a class $x$ by $\Delta _i x$. The
natural operators $\Delta_i$ are called \emph{higher $BV$
operators}.
\end{defin}

The action of $g\in H_*(SO(n))$ on $b_0 \in H_0(S^{n-1})$ is
denoted by  $g\cdot b_0$. In particular, for $n=2$, $\Delta _1$ is
the Batalin-Vilkovisky operator on the homology of $2$-fold loops
on $SO(2)$-spaces \cite{Getzler:BVAlg}, simply denoted by $\Delta$
or $BV$.

\sm
For $n>0$, consider the natural (unpointed) action of $SO(n)$ on $S^{n-1}$. For $g\in H_* SO(n)$ and $x\in H_* S^n$, we denote by $g.x$ the image of $g\otimes x \in H_* S^{n-1}$ by the map induced in homology by the natural action. Recall that $b_0\in H_0 S^{n-1}$ is a generator. The element $g.b_0$ is none but $H_*(ev)(g)$ where $ev:SO(n) \lra S^{n-1}$ is the evaluation of the action at the point $(0,\cdots,0,1)$ of $S^{n-1}$. This can classically be computed as the edge
homomorphism in the Serre spectral sequence of the fibration
sequence $SO(n-1)\hookrightarrow SO(n) \lra S^{n-1}$. We get:
\begin{eqnarray}
\label{sonsursn} e_{n-1}\cdot b_0=b_{n-1} ~~~~
\end{eqnarray}
As indicated in \cite[p. 291-299]{hatch}, $SO (n)$ possesses a
nice cell structure compatible with the reflection map $\kappa$. The
{\it modulo} $2$ cellular chain complex ${C_*^{\mathrm{cell}}(SO
(n))}$ is the differential graded bialgebra $\Lambda(d^1,d^2,
\ldots, d^{n-1})$ with zero differential and Cartan diagonal on
the generators $d^i$.

\sm
There is a chain equivalence
$$
C_*^{\mathrm{cell}}(SO(n)\times S^{n-1})\cong C_*^{\mathrm{cell}}(SO(n))\otimes C_*^{\mathrm{cell}}(S^{n-1})
$$
and we will need a cellular approximation of the action
$SO(n)\times S^{n-1}\to S^{n-1}$ and its effect
$$
C_*^{\mathrm{cell}}(SO(n))\otimes C_*^{\mathrm{cell}}(S^{n-1}) \to C_*^{\mathrm{cell}}(S^{n-1}) ~~.
$$
on cellular chains.
\begin{lem}
\label{lemlemlem} There is a cellular approximation of the action
$SO(n)\times S^{n-1}\to S^{n-1}$ which induces on cellular chains
the map
$$
C_*^{\mathrm{cell}}(SO(n))\otimes C_*^{\mathrm{cell}}(S^{n-1}) \to C_*^{\mathrm{cell}}(S^{n-1}) ~~.
$$
given by
\begin{itemize}
\item $d^i. e^j=0$ except if $j=0$ or $i=0$,
\item $d^{0}.e^i = e^i$ for $i\geq 0$,
\item $d^i. e^j=0$ for $i>0$,
\item $d^i. Te^j = T(d^i. e^j)$ for $0\leq i\leq n-1$ and $ 0\leq j \leq n-1$.
\end{itemize}
where the action of the generator $d^i\in C_*^{\mathrm{cell}}SO
(n)$ on the generator $e^j \in C_*^{\mathrm{cell}}S^{n-1}$ is
denoted by $d^i. e^j$.
\end{lem}

\sm NB: $d^i$ is a generator in the chain complex
$C_*^{\mathrm{cell}} SO(n)$, while $d_i$ is the induced homology
class in $H_* SO(n)$.


%
{\bf Proof.~}
We prove the lemma by induction on $n\geq 2$. We begin by the
construction of a nice cellular approximation of the action for
each $n$.

\sm
We first observe that there is a commutative diagram for $k\geq n$
\begin{eqnarray}
\label{champions}
\xymatrix{
SO(n)\times S^{n-1}\ar[r]\ar[d]& S^{n-1}\ar[d]\\
SO(k)\times S^{k-1}\ar[r]& S^{k-1}
}
\end{eqnarray}
where the vertical maps are the usual inclusions, which are
cellular.


\sm Beginning with $n=2$,  where $SO(n)= S^1$ and the action of
$SO(n)$ on $S^1$ is the action by translation of $S^1$ on itself,
we construct inductively for $n\geq 2$ compatible cellular
approximations in the sense that Diagram (\ref{champions})
commutes, and such that the restriction of the action $SO(n)\times
S^{n-1}\to S^{n-1}$ to $\{d^0\}\times S^{n-1}$ is the identity.
This follows by a simple application of the cellular approximation theorem for CW-pairs.

\sm
As our cellular
approximations restrict to the identity on $SO(n)\times S^{n-1}\to
S^{n-1}$, we immediately derive that
$$
d^{0}.e^i = e^i
$$
for $0\leq i\leq n-1$ and all $n\geq 2$ at the cellular chains
level.

\sm Now assume that  Lemma \ref{lemlemlem} holds for some $n\geq 2$, we want to derive that it holds for $(n+1)$ as well.
As we have chosen compatible cellular approximations of the action,
we have for $i, j\leq n-1$ that $d^i.e^j=0$ if neither $i\neq 0$ nor $j\neq 0$. Moreover
$$
d^{i}.e^0=e^i+Te^i~,~~ i>0;  ~~~~~~~~ d^{0}.e^i = e^i
$$
for $i\leq n-1$.
By degree reasons we have for $i, j>0$
$$
d^i.e^n =  d^n. e^j =0
$$
It remains thus to check:
$$
d^{n}.e^0=e^n+Te^n.
$$
But by Formula (\ref{sonsursn}), we must have $d^{n}.e^0=e^n+Te^n$
because $e^n+Te^n$ is the unique $n$-cycle that represents the top
class of $H^* S^n$.

\sm
The last point of Lemma \ref{lemlemlem} follows from similar considerations.

\sm We still have to check the case $n=2$ as  initial input of the induction process.
As already noticed, in this
case $SO(n)= S^1$, and the action of $SO(n)$ on $S^1$ is the
action by translation of $S^1$ on itself. The effect of the action
in homology is simply given by the algebra structure of $H^* S^n =
\Lambda (b_1)= \Lambda (e^1)$, and we have $d_1.b_0= b_1$ in
homology, which forces at the cellular chains level $d^1.e^0=
e^1+Te^1$ (but this is actually obvious in this case at the
topological level). This settles the case $n=2$.

\subsection{The fundamental diagram}

Let $Z$ be an $\fdn$-space, {\it i.e.} a space with an action of the
operad $\fdn$ of framed little $n$-disks. Then  $Z$ is in
particular a pointed $SO(n)$-space with an action of the little
$n$-disks operad, such that the following diagram commutes:
\begin{eqnarray}
\label{facto1}
\xymatrix{
SO(n)\times S^{n-1}\times Z \times Z \ar[rrrr]^{{\act}^\mathrm{D}}\ar[dd]^{id \times \theta}&&&& S^{n-1} \times Z\times Z \ar[dd]^{\theta}\\
&&&&\\
SO(n)\times Z \ar[rrrr]^{{\act}}&&&& Z }
\end{eqnarray}
Here $\theta$ is the action of little $n$-disks, $\act$ is the
action of $SO(n)$ on $Z$, and ${\act}^\mathrm{D}$ is the diagonal
action of $SO(n)$, where the action of $SO(n)$ on $S^{n-1}$ is the
natural one (i.e. the restriction of the action of the natural
action of $SO(n)$ on $\mathbb{R}^n$ to its unit sphere $S^{n-1}$).
We emphasize again that we can use, instead of the space of $2$
small disjoint little $n$-disks in a big $n$-disk, the
equivariantly homotopy equivalent space $S^{n-1}$.

\sm We have to translate the commutation of Diagram (\ref{facto1})
in terms of relations between operations coming from the $H_* SO
(n)$ action (the $\Delta_i$'s) and the operations from the loop
structure (Pontryagin product, Browder bracket, Kudo-Araki
operations). The relations involving Kudo-Araki operations are (of
course) a little more difficult to derive. We therefore postpone
this calculation to the Section \ref{bv:dl}.

\subsection{Higher $BV$ operations and Pontryagin products}
\label{bv:pontryagin}

We now compute the behavior of the $H_*SO (n)$ action on $H_*
Z$ with respect to Pontryagin products in $H^*Z$. The computation is:
{\small
\begin{eqnarray}
\Delta_i (x*y)&=& {\act}_*(d_i \otimes \theta_* (e_0 \otimes x \otimes y))\\
&=&{\act}_*(\mathrm{Id}_{H_*(SO (n))}\otimes \theta_* ) (d_i\otimes(e_0\otimes x\otimes y))\nonumber\\
&=&\theta _*(\sum_{\alpha+\epsilon+\gamma=i}d_\alpha . e_0\otimes d_\epsilon .x\otimes d_\gamma.y)\nonumber\\
&=& \delta^i_{n-1}\theta _*(e_{n-1}\otimes x \otimes y)+ \theta_*(\sum_{\epsilon+\gamma=i}e_0\otimes d_{\epsilon}.x \otimes d_\gamma y) \nonumber\\
&=&\delta^{i}_{n-1}\{x, y\} +\theta_*( e_0\otimes  \sum_{\epsilon+\gamma=i} d_{\epsilon}.x\otimes d_{\gamma}. y)\nonumber \\
&=&\delta^{i}_{n-1}\{x, y\}
+\sum_{\epsilon+\gamma=i}(\Delta_\epsilon x *\Delta_\gamma
y)\nonumber
\end{eqnarray}}
Hence we recover (see \cite{Salvatore-Wahl:FrameddoBVa}):
\begin{theor} \label{theo:bat} Let $Z$ be an
$f\mathcal{D}_n$-algebra and $x,y$ be elements in the modulo $2$
cohomology of $Z$. Then, for $i\leq n-1$:
\begin{eqnarray}
\Delta_i (x*y)&=&\delta^{i}_{n-1}\{x, y\}
+\sum_{\epsilon+\gamma=i}\Delta_{\epsilon}x*\Delta_{\gamma} y
\quad .
\end{eqnarray}
\end{theor}
We notice that $H_*SO(n)$ has a primitive element $p_i$ in every
odd degree $i$. Moreover, $p_i$ is equal to $d_i$ up to
decomposable elements.  For $i$ odd, let $\Delta_{p_i}$ be the
operator associated to the primitive class $p_i$. For example
$$
p_1= d_1, ~~ p_3= d_3+d_2*d_1
$$
hence
$$
\Delta_{p_1}= \Delta_1,~~ \Delta_{p_3}= \Delta_3+\Delta_2
\Delta_1=\Delta_3+\Delta_1
\Delta_2 \quad .
$$
Using the same method as in the proof of Theorem
\ref{theo:bat}, we see that for $i$ odd with $i< n-1$, then
$\Delta_{p_i}$ is a derivation with respect to the Pontryagin
product, while for $n$ even, $\Delta_{p_{n-1}}$ is a $BV_n$
operator inducing the natural Gerstenhaber algebra structure of
the cohomology.
\begin{theor}
\label{theo:bat:1} Let $Z$ be an $\fdn$-algebra and $x,y$ be
elements in the modulo $2$ cohomology of $Z$. For any odd integer $i$ with $i\leq
n-1$:
\begin{eqnarray}
\Delta_{p_i} (x*y)&=&\delta^{i}_{n-1}\{x, y\} +(\Delta_{p_i}x)y +x
(\Delta_{p_i} y) \quad .
\end{eqnarray}
\end{theor}
For $i=n-1$ (that is, for $n$ is even), this is the {\it modulo}
$2$ expression of the relations defining a $BV_n$-algebra (see
Definition \ref{definition BV_n-algebre}). We get:
\begin{cor}
Let $Z$ be an $\fdn$-algebra. For $i\leq n-1$, $i$ odd, let
$\Delta _{p_i}$ be the natural operator on $H_* Z$ obtained as the
action of the primitive element $p_i \in H_* SO$ on $H_* Z$. Then
for $i<n-1$, the $\Delta_{p_{i}}$ operator is a derivation with
respect to the Pontryagin product. If $n$ is even, $\Delta_{p
_{n-1}}$ endows $H_* Z$ with the structure of a $BV_n$-algebra.
\end{cor}

\subsection{Higher BV operations and Browder brackets}
\label{bv:browder}

We now compute the behavior of the $H_*SO (n)$ action on $H_* Z$ with respect to Browder brackets. The computation is
essentially the same as above. Let $x, y \in H_* Z$ {\small
\begin{eqnarray}
\Delta_i (\{x,y\})&=& {\act}_*(d_i \otimes \theta_* (e_{n-1} \otimes x \otimes y))\\
&=&(\mathrm{Id}_{H_*SO (n)}\otimes \theta_*) (d_i\otimes e_{n-1}\otimes x\otimes y)\nonumber\\
&=& \theta_*(\sum_{\alpha+\epsilon+\gamma=i}d_\alpha.e_{n-1}\otimes d_\epsilon .x\otimes d_\gamma.y)\nonumber\\
&=& \sum_{\epsilon+\gamma=i}e_{n-1}\otimes d_{\epsilon}.x \otimes
d_\gamma y \nonumber
\end{eqnarray}}
Hence we get:
\begin{theor}
\label{theo:bat:2} Let $\Omega^n X$ be an $n$-fold loop space and
$x,y$ be an elements in the modulo $2$ cohomology of $\Omega^n X$.
For $i\leq n-1$:
\begin{eqnarray}
\Delta_i (\{x,y\})&=&\sum_{\alpha+\epsilon=i}\{\Delta_\alpha x,
\Delta_\epsilon y \} \quad .
\end{eqnarray}
\end{theor}

\section{Higher $BV$ operators and Kudo-Araki operations}
\label{bv:dl}

Let $n\geq 2$. We now compute the relation between $\Delta_i$ and
$Q_j$ operations, for $i,j\leq n-1$. Let $Z$ be any algebra over
the framed $n$-disks operad $\fdn$. Typically, $Z=\Omega^n X$ is
an $n$-fold loop space.

\subsection{Decomposition of the $ SO(n) $-action}
The topological group $ SO(n) $ acts on
$S^{n-1}$ in the standard fashion,  and acts diagonally on
products of $SO(n)$-spaces. Thus we obtain a natural action of $
SO(n) $ on $S^{n-1}\times Z\times Z$, and this action is denoted
by ${\act}^\mathrm{D}$. We break down this action into two
actions of $ SO(n) $, the first one is an action $a_1$ on $S^{n-1}$, and one by
$a_2$ on $Z\times Z$. In other words, the diagonal action
${\act}^\mathrm{D}$ of $\sodeux$ decomposes as:
{\footnotesize
\begin{eqnarray}
~~~~~~\label{dec:act:diag} ~~~~~~~~SO(n) \times S^{n-1} \times Z^2
\stackrel{\mathrm{D}\times Id}{\lra} SO(n)^2 \times S^{n-1} \times
Z^2 \stackrel{1\times a_1}{\lra}  SO(n) \times S^{n-1} \times Z^2
\stackrel{a_2}{\lra} S^{n-1} \times Z \times Z~~~~~~
\end{eqnarray}
}

\sm  These two actions $a_1$ and $a_2$ are $\zdeux$-equivariant.
That is,  Diagram (\ref{dec:act:diag}) is a diagram of
$\zmoddeux$-spaces. Hence the two actions $a_1$ and $a_2$ pass to
the quotient $\etendnnz$, on which  they induce corresponding
actions (denoted by $\bar{a}_1$ and $\bar{a}_2$). We therefore get
a commutative diagram: {\tiny
\begin{eqnarray} \label{facto2}
\xymatrix{
SO(n) \times S^{n-1} \times Z^2\ar[rrr]^{{\act}^{\mathrm{D}}}\ar[rd]^{\mathrm{D}\times Id}\ar[dd]^{Id \times \zeta}&&&SO(n) \times S^{n-1} \times Z^2 \ar[r]^{~~~~~~~\theta}\ar[dd]^{\zeta}&Z\\
& SO(n)^2 \times S^{n-1} \times Z^2\ar[r]^{1\times a_1}\ar[d]^{Id \times \zeta}&  SO(n) \times S^{n-1} \times Z^2\ar[d]^{Id \times \zeta}\ar[ru]^{a_2}& \\
SO(n) \times \etendnnz\ar[r]^{\mathrm{D}\times Id} &SO(n)^2\times \etendnnz
\ar[r]^{1\times \bar{a}_1}& SO(n) \times \etendnnz  \ar[r]^{
\bar{a}_2} &\etendnnz\ar[ruu]_{\bar{\theta}}& }
\end{eqnarray}
}We have a decomposition of the diagonal action ${\act}^\mathrm{D}$ as:
$$
SO(n)  \times \etendnnz \lra  SO(n) ^2 \times \etendnnz
\stackrel {\relax 1\times \bar{a}_1} {\lra}  SO(n)  \times \etendnnz
\stackrel{\bar{a}_2}{\lra}\etendnnz
$$
By the commutation of Diagram (\ref{facto1}) and  Diagram
(\ref{facto2}), we only need to compute
\begin{eqnarray}
\label{equationcalcul1} \Delta_i Q_j x =
\bar{\theta}_*(\bar{a}_2)_*(1\otimes (\bar{a}_1)_*)(1\otimes
\mathrm{D}_*)(d_i \otimes [e^j \otimes x \otimes x ])
\end{eqnarray}
hence the main issue is to compute $(\bar{a}_2)_*$ and
$(\bar{a}_1)_*$.

\subsection{\relax Calculation of $(\bar{a}_1)_*$}

The result is as follows:
\begin{proposition}
\label{calcula1} The following formul{\ae} hold:
\begin{itemize}
\item For  $x\neq y$,
$$
(\bar{a}_1)_* (d_i \otimes [e^j \otimes x\otimes y]) = [d_i*e^j
\otimes x\otimes y]
$$
where $d_i*e^j=0$ except in the cases $d_0*e^i= e^i$, $d_j *e^0= e^j$ and
$d_{n-1}*e^0= e^{n-1}$.
\item For $x=y$,
$$
(\bar{a}_1)_* (d_i \otimes [e^j \otimes x\otimes x]) = [d_i*e^j
\otimes x\otimes x]
$$
with $d_i*e^j=0$ except in the case $d_0*e^i= e^i$.
\end{itemize}
\end{proposition}

{\bf Proof.~}We first notice that the formula for $(\bar{a}_1)_*(d_{i} \otimes
[e^0\otimes x \otimes y])$ holds because such classes come from
$H_*  (S^{n-1}\times Z\times Z)$, and we know how the action of
$H_* SO (n)$ is  at this level (Formula (\ref{sonsursn})), as well
as the map $H_* (S^{n-1}\times Z\times Z\lra \etendnnz)$ (Formula
(\ref{quotientmapcomput})). In the same way the formula for
$(\bar{a}_1)_*(d_{i} \otimes[e_{n-1}\otimes x \otimes y])$ for
$x\neq y$ must hold.

\sm In general, we can not argue this way to compute $
(\bar{a}_1)_* (d_i \otimes [e^j \otimes x\otimes x])$ for $0<i<n$.
Instead, we rely on the naturality of the isomorphisms
(\ref{isos2}). Let $\tilde{x}$ be a cycle representing $x$, and
recall $d^i$ is a generator in the cellular complex of $SO(n)$
corresponding to the cohomology class $d_i$. The class $d_i
\otimes [e^j \otimes x\otimes x]$ is represented by the cycle $d^i
\otimes e^j \otimes \tilde{x} \otimes \tilde{x}$. At the level of
chains, $(a_1)_* (d^i \otimes e^j \otimes \tilde{x} \otimes
\tilde{x})= d^i.e^j \otimes \tilde{x} \otimes \tilde{x}$, which by
the Lemma \ref{lemlemlem} is zero except  if $j=0$. In this case,
we compute
$$
(a_1)_* (d^i \otimes e^0\otimes \tilde{x} \otimes \tilde{x})= d^i.e^0 \otimes \tilde{x} \otimes \tilde{x}= (e^i+Te^i)\otimes \tilde{x} \otimes \tilde{x}
$$
But now
$$
[(e^i+Te^i)\otimes x \otimes x]=[e^i\otimes x \otimes x]+  [Te^i\otimes x \otimes x]= [e^i\otimes x \otimes x]+[e^i\otimes x \otimes x]=0
$$
in $H_*(\etendnnz)$.
\hfill $\square$

\subsection{Calculation of $(\bar{a}_2)_*$}

We can decompose $a_2$ a little further. The map $\bar{a}_2$ is
the bottom composition in the commutative diagram:
$$
\xymatrix{
SO(n)\times S^{n-1} \times Z^2 \ar[r]\ar[d]& S^{n-1} \times (SO (n) \times Z)^2 \ar[r] \ar[d]& S^{n-1} \times Z^2 \ar[d] \\
SO(n) \times (S^{n-1} \times_{\zdeux} Z^2)\ar[r]_{\varphi} & S^{n-1}
\times_{\zdeux} (SO(n)\times Z)^2 \ar[r]_{~~~\psi}& S^{n-1} \times
_{\zdeux} Z^2 }
$$
We can compute $\psi_*$ easily:
\begin{eqnarray}
\label{psistar} \psi([e^j \otimes (s_1\otimes  x)\otimes
(s_2\otimes y)]) = [e^j \otimes s_1.x\otimes s_2.y]~~.
\end{eqnarray}
This relies as above, on the naturality in $H_*Z$ of the Serre
spectral of the fiber  sequence:
$$
Z\times Z \lra \etendnnz \lra \rp ^{n-1}
$$
at the $E_2$-term, together with the naturality explained in
Section \ref{natur}.


\subsection{Calculation of $\varphi _*$}

\noindent\subsubsection*{The result.} The result for the computation of $\varphi _*$ is as follows.

\begin{lem}
For all $0\leq r\leq {n-1}$, $\alpha \in H^* SO(n)$, $x\in
H^* Z$
\begin{eqnarray}
\label{formulegenerale} \varphi_* (\alpha\otimes  [e^r
\otimes x \otimes x]) &=&  \sum_{i, r+2i-|\alpha|\leq {n-1}}
[e^{r+2i-|\alpha|}\otimes (Sq^i_* \alpha \otimes x) \otimes
(Sq^i_* \alpha \otimes x)] \\
\nonumber &&+\delta_{n-1}^r\sum_{|\alpha'|<|\alpha"|}[e^{n-1}
\otimes (\alpha' \otimes y) \otimes (\alpha"\otimes y)]
\end{eqnarray}
where $D_* (\alpha)= \sum \alpha'\otimes \alpha"$ and $D$ is the
diagonal of $SO(n)$.
\end{lem}
We will proceed as follows: the formula up to a defect term
appears already elsewhere as we shall briefly recall;  we then
compute this defect term using the map $\xi$ of Section \ref{recoll:dyer:lashoff}.

\noindent\subsubsection*{Related results in the literature and the
method of calculation.} Our method to compute $\varphi_* $ is to a
certain extent classical. Very close computations occur for
instance in \cite[p. 363, Theorem 3.2.
(ii)]{Cohen-Lada-May:homiterloopspaces} or \cite[p.
11-13]{Cohen-Lada-May:homiterloopspaces}.
We consider the diagonal map
$$
\varphi : X\times (S^{n-1} \times_{\zdeux} Y^2)\to
S^{n-1} \times_{\zdeux}(X \times Y)^2
$$
induced by the natural map
$$
(\shuffle)\circ (D \times Id_{S^{n-1}\times Y^2 }) : X\times
S^{n-1} \times Y^2\to  S^{n-1} \times (X \times Y)^2
$$
where $D$ is the diagonal of $X$ and ${\shuffle}$ is the obvious
shuffling of factors in the product $X\times X \times
S^{n-1}\times Y \times Y$. We are interested in the special case
$X= SO(n)$. We will actually compute the map $ \varphi_* : H_*(Y)
\otimes H_* \etennd X \lra H_* \etennd (X\times Y)$ in general.
According to  \cite[Lemma 4.2 (i), p. 368, Theorem
3.2]{Cohen-Lada-May:homiterloopspaces}:
\begin{lem} There is a
formula, for $x\in H_* X$ and $y\in H_* Y$:
\begin{eqnarray}
\label{formulegenerale:cohen} \varphi_* (x\otimes  [e^r \otimes y
\otimes y]) &=&  \sum_{i, r+2i-|x|\leq n-1}
[e^{r+2i-|\alpha|}\otimes (Sq^i_* x \otimes y) \otimes (Sq^i_* x
\otimes y)] + \delta_r^{n-1}\Gamma_r
\end{eqnarray}
where $\Gamma_r$ is in the image of map induced in homology by the
quotient map $S^{n-1}\times (X\times Y)^2 \to S^{{n-1}}\times_{\zdeux}
(X\times Y)^2 $
\end{lem}
We warn the reader that our notations  differ from those of the
reference, and moreover that the formula \cite[Lemma 4.2, (i), p.
368, Theorem 3.2]{Cohen-Lada-May:homiterloopspaces} has a small
typo in the indices. The fact that $\Gamma_r$ is zero for $r< n-1$ is not part of  \cite[Lemma 4.2 (i), p. 368, Theorem
3.2]{Cohen-Lada-May:homiterloopspaces}. This can be seen as follows. We can apply homology to the commutative diagram:
{\tiny
$$
\xymatrix{
X \times \etenrd
(Y)\ar[rr]^{\varphi}\ar[dd]^{id\times \xi}\ar[rd]^{\iota}&& \etenrd  (X\times
Y)\ar@{-}[d]\ar[rd]^{\iota}&\\
&X \times \etennd (Y)\ar[rr]^{~~~~~~~\varphi}\ar[dd]\ar'[d][dd]^{id\times \xi}&\ar[d]^{\xi}& \etennd  (X\times
Y)\ar[dd]\ar'[d][dd]^{\xi}\\
X\times  (S^{r} \wedge  Y\times  Y)\ar@{-}[r]\ar[rd]^{\iota}&\ar[r]^{\tilde{D}}& S^{r}
\wedge  [(X\times  Y)\times  (X\times  Y)]\ar[rd]^{\iota}&\\
&X\times  (S^{n-1} \wedge  Y\times  Y)\ar[rr]^{\tilde{D}~~~}&& S^{n-1}
\wedge  [(X\times  Y)\times  (X\times  Y)]&,
}
$$
}
where the maps labelled $\iota$ are induced by natural inclusions.
The element $x\otimes  [e^r \otimes y
\otimes y]$ comes from the homology of $X \times \etenrd (Y)$
hence the element $\Gamma_r$ comes from some element $\Gamma'_r$ in the homology of $\etenrd  (X\times
Y)$. Brackets are detected by the map $\xi$ (Formulas \ref{sacre:formules}), in the sense that $\xi_*$ is injective on the subspace generated by bracket elements. We see that
$\Gamma_r$ is non zero if and only if the image of $\Gamma'_r$ by the map
$$
(\iota \xi)_*: H_* \etenrd  (X\times
Y) \to H_* S^{r}
\wedge [ (X\times  Y)\times  (X\times  Y)]\to  S^{n-1}
\wedge  [(X\times  Y)\times  (X\times  Y)]
$$
is non zero.  But the map $\iota_*
$ is zero on the image of $\xi_*
$. This means that $\Gamma'_r$ is zero for $r<n-1$, and
for $r=n-1$:
$$
\Gamma:= \Gamma_{n-1} =\sum [e^{n-1} \otimes x_1 \otimes y_1 \otimes x_2
\otimes y_2 ]
$$
We set
$$
G= \sum_{i, r+2i-|x|\leq {n-1}} [e^{r+2i- |x|}\otimes (Sq^i_*
x\otimes y) \otimes (Sq^i_* x \otimes y)]
$$
So that $\varphi_* (x\otimes  [e^{n-1} \otimes y \otimes
y])=G+\Gamma$.

\noindent\subsubsection*{Computation of $\Gamma$.} There is a
commutative diagram
$$\xymatrix{ X \times \etennd
(Y)\ar[r]^{\varphi}\ar[d]^{id\times \xi}& \etennd  (X\times
Y)\ar[d]^{\xi}\\
X\times  [S^{n-1} \wedge (Y\times  Y)]\ar[r]^{\tilde{D}}& S^{n-1}
\wedge [ (X\times  Y)\times  (X\times  Y)] }
$$
Here $\xi$ is as in Section \ref{otherview} and $\tilde D$ is
induced by the diagonal on $X$ and shuffling the factors.
We compute that
\begin{eqnarray}
\xi_* \varphi_*  (x\otimes  [e^{n-1} \otimes y \otimes y])&=&{\xi}_*
(G+\Gamma)\\
\nonumber &=&(\tilde{D})_*(id\otimes \xi_*)(x\otimes  [e^{n-1}
\otimes y \otimes y])
\end{eqnarray}But
\begin{eqnarray}
(\tilde{D})_*(id\otimes \xi_*)(x\otimes  [e^{n-1} \otimes y
\otimes y])&=&(\tilde{D})_*(x\otimes  b_{n-1} \otimes y \otimes y)\\
\nonumber &=&\sum b_{n-1} \otimes x' \otimes y \otimes x" \otimes
y
\end{eqnarray}where $D_*  x = \sum x'\otimes x"$. If we set $D_*  x = \sum
x'\otimes x" = \sum_{x'= x"}x'\otimes x" + \sum_{x'\neq
x"}x'\otimes x"$ then if $|x|$ is even
$$
Sq^{|x|/2}_* (x )= \sum_{x'= x"}x'\otimes x"
$$
because $Sq^{|x|/2}_*$ is dual to the cup square and the diagonal $D_*$
is dual to the cup product. In other words
{\footnotesize
\begin{eqnarray}
(\tilde{D})_*(id\otimes \xi_*)(x\otimes  [e^{n-1} \otimes y
\otimes y])&=& \sum_{x'= x"}e_{n-1} \otimes (Sq^{|x|/2}_* x \otimes
y) \otimes (Sq^{|x|/2}_* x
 \otimes  y) \\
&&+ \sum_{x'\neq x"} e_{n-1}
\otimes x' \otimes y \otimes x" \otimes  y \nonumber
\end{eqnarray}
}
where the first term has to considered as zero if $x$ has odd degree.
We compute ${\xi}_* (G)$ is non zero if and only if $x$ has
even degree (Formula (\ref{sacre:formules})), in which case we
have
\begin{eqnarray}
{\xi}_* (G)&=&{\xi}_* \sum_{i, r+2i-|x|\leq {n-1}} [e^{r+2i-
|x|}\otimes (Sq^i_* x\otimes y) \otimes (Sq^i_* x \otimes
y)]\\
\nonumber &=&{\xi}_*[e^{n-1}\otimes (Sq^{|x|/2}_* x\otimes y)
\otimes
(Sq^{|x|/2}_* x \otimes y)]\\
\nonumber &=&b_{n-1}\otimes (Sq^{|x|/2}_* x\otimes y) \otimes
(Sq^{|x|/2}_* x \otimes y)
\end{eqnarray}

\sm So it follows that
$$ {\xi}_* (\Gamma)=
\sum_{x'\neq x"} b_{n-1} \otimes x' \otimes y \otimes x" \otimes y
$$
In view of the formulas for $\xi_*$ (\ref{sacre:formules}), we
must have
$$
\Gamma =\sum_{|x'|< |x"|} [e^{n-1} \otimes x' \otimes y \otimes
x" \otimes y].
$$

\subsection{Higher BV Relations}

Let $j>0$. Putting together the previous calculations, we get:
\begin{eqnarray}
\label{calculfinal} \Delta_i Q_{j}x &= &
\bar{\theta}_*(\bar{a}_2)_*(1\otimes (\bar{a}_1)_*)(
\mathrm{D}_*\otimes 1)(d_i \otimes [e^j \otimes x \otimes x ])\\
\nonumber&=&\bar{\theta}_*(\bar{a}_2)_*(1\otimes (\bar{a}_1)_*)
(\sum_{m+l=i} (d_m \otimes d_l )\otimes [e^j \otimes
x\otimes x ])\\
\nonumber&=&\bar{\theta}_*(\bar{a}_2)_*(\sum_{m+l=i} d_m \otimes
[d_l.e^j \otimes
x\otimes x ])\\
\nonumber&=&\bar{\theta}_*(\bar{a}_2)_* (d_i \otimes [e^j \otimes
x\otimes x ])\\
\nonumber&=&\bar{\theta}_*(\psi_*\varphi_*) (d_i \otimes [e^j
\otimes x\otimes x ])\\
\nonumber&=&\bar{\theta}_*\psi_*(\sum_{k, j+2k-i \leq n}
[e^{j+2k-i}  \otimes (Sq^k_* d_i \otimes x)\otimes(Sq^k_* d_i
\otimes x)] +\\
\nonumber&& \delta_j^{n-1}[e_{n-1}\otimes \sum_{\alpha+\beta=i,
\alpha<\beta}  (d_\alpha \otimes  x) \otimes  (d_\beta \otimes
x)]~)
\end{eqnarray}
Recall that  ${\act}$ is the action of $SO(n)$ on $Z$, and by $\Delta_ix$ the element ${\act}_* (d_i \otimes x)$.
If we let $Sq^k_* \Delta_i x$ denote ${\act}_* (Sq^k_* d_i
\otimes x)$, we get:
\begin{eqnarray}
\label{prehigherbv} \Delta_i Q_{j}x= \sum_{k, j+2k-|x| \leq j}
Q_{j+2k-i} (Sq^k_* \Delta_i
x)+\delta_j^{n-1}\sum_{\alpha+\beta=i, \alpha<\beta}
\{\Delta_\alpha x, \Delta_\beta x\}
\end{eqnarray}
The action of the Steenrod algebra in our case is well known:
$$
Sq^k_* d_i = C^{i-k}_k d_{i-k}~~ \mathrm{for~}2k\geq i~~.
$$
This comes from the action of the Steenrod algebra on $H_*\rp ^n$. Here $C^{i-k}_k$ denotes the usual binomial coefficient
$k!/(i-k)!i!$ reduced modulo $2$.

\begin{theor}
The following formula holds in the homology of $\fdn$-spaces:
\begin{eqnarray}
\label{higherbv} \Delta_i Q_{j}x= \sum_{k, j+2k-i \leq n, 2k \geq i}
C^{i-k}_k Q_{j+2k-i} (\Delta_{i-k}
x)+\delta_j^{n-1}\sum_{\alpha+\beta=i, \alpha<\beta}
\{\Delta_\alpha x, \Delta_\beta x\}
\end{eqnarray}
\end{theor}
In particular, if $\Delta_i$ does not support any Steenrod
operations, for example if we consider the operation associated to
the primitive $p_{2^i-1} \in H_* SO(n)$, we get:
$$
\Delta_{p_{2^i-1}}Q_j x = Q_{j-2^i+1} \Delta_{p_{2^i-1}} x+
\delta_j^{n-1} \{x, \Delta_{p_{2^i-1}}x \}
$$
We insist again  that these formulas hold for any $\fdn$-algebra.
By specializing the result to the case $n=2$, $i=j=1$, we obtain:
\begin{cor}
\label{comm:bv:q} Let $\Omega^2 X$ be a $2$-fold loop space, and
$x$ be an element in  $H_*\Omega^2 X$
\begin{eqnarray}
\label{bvformula} \Delta (Q_1 x)= \{\Delta (x), x\}+ \Delta (x)* \Delta
(x)= \Delta (x*(\Delta x))~~~.
\end{eqnarray}
\end{cor}
Where $\Delta := \Delta _1 = BV $ is the Batalin-Vilkovisky operator. Here we have considered $X$ as a trivial $SO(2)$-space.
We insist again that such a formula holds in the {\it modulo}
$2$ homology of any $\fd _2$ algebra.

In particular if $Y$ is  an infinite loop space, all the Browder
brackets on $H_*Y$ are trivial, all Kudo-Araki operations are
defined, and there is an action of the infinite special orthogonal
group $SO$. For example, in the case of infinite loop spaces, we have:
\begin{eqnarray}
\Delta_i Q_j x = Q_j \Delta _{i/2} x \quad ~\mathrm{for~} 2i>j.
\end{eqnarray}
Here $\Delta _{i/2}$ is zero if $i$ is odd.
For example $\Delta _{2p+1} Q_0$ is trivial. In other words, the odd
generators of $H_* SO$ act trivially on the Pontryagin squares in
the homology of infinite loop spaces. More generally, as
$H_*SO(n)= \Lambda (d_1, \ldots , d_{n-1})$, any homogeneous
element of odd degree is a sum of homogeneous monomials, each of
them having some $d_{i}$ factor with $i$ odd. Hence, each such
monomial act trivially on Pontryagin squares.

\begin{cor}
\label{cor:nul:sur:square} Let $Y$ be an infinite loop space. The
action of $H_{\mbox{odd}} SO$ on $H_* Y$ is trivial on Pontryagin
squares.
\end{cor}


\begin{rem}
\label{remmay} We have noticed in previous work that the action of
$SO(n)$ on an $n$-fold loop space can be interpreted  in terms of
$\circ$-products and of the $J$-homomorphism \cite[section
5]{gaudensmeni0}. Formulas relating the circle product with
elements coming from $H_* Q_1 S^0$ (whom the $J$-homomorphism
factors through) and Kudo-Araki operations coming from the loop
structure are already known (see for example Formula
(\ref{actionhj})). These formulas come from the structure of
$E_\infty$-module over the $E_\infty$-ring $QS^0$. Albeit less
general (as this works only for elements coming from the
$J$-homomorphism), the new set of formulas given above is new and
quite simple. It relies on the structure of
$f\mathcal{D}_n$-space. It should be in principle possible to get
our formulas from   Formula (\ref{actionhj}) directly, by using an
explicit description of the generators of the  homology image of
$J$. This looks however quite tedious. We will go a little further
in this direction in Section \ref{comp:hur:qso}.
\end{rem}


\section{Calculations of the Batalin-Vilkovisky structures of $H_*(\Omega^2 S^k)$ for $k\geq 2$}
\label{sec:bvspheres}

In this section, we accomplish the promised computations. As before, the notation $*$ stands for Pontryagin products.

\subsection{Classical splittings}
We denote by $\eta: S^3\longrightarrow S^2$, $\nu: S^7
\longrightarrow S^4$, and $\sigma: S^{15}\longrightarrow S^8$ the
classical maps of Hopf invariant $1$. Let $a\in \pi_1 SO(2)$,
$b\in \pi_3 SO(4)$ and $c\in \pi_7 SO(8)$ be generators of the
$\mathbb{Z}$-summands that map to the Hopf maps in the
corresponding degrees under the $J$-homomorphism. We denote the
corresponding stable elements in the homotopy of $SO$ and $QS^0$
in the same way.
\begin{lem}\label{splitting lemma}
Let $F\buildrel{j}\over\hookrightarrow
E\buildrel{p}\over\twoheadrightarrow B$ be a homotopy fiber
sequence with $F$, $E$ and $B$ path-connected. Suppose that there
is a map $s:B\rightarrow E$ such that the composite $p\circ s$ is
a homotopy equivalence and that $E$ is an H-space with
multiplication $\mu$. Then the composite $\mu\circ (j\times
s):F\times B\buildrel{\simeq}\over\rightarrow E$ is a homotopy
equivalence. In particular $j$ admits a retract up to homotopy.
\end{lem}
Let $f=\eta$, $\nu$ or $\sigma:S^{2d-1}\rightarrow S^d$,
$d\in\{2,4,8\}$, be a Hopf fibration. Let $\partial:\Omega
S^d\rightarrow S^{d-1}$ be the connecting homomorphism. Let
$\iota:S^{d-1}\longrightarrow  \Omega S^d$ be the adjoint of the
identity map of $S^d$. Since $\pi_{d-1}(\partial)$ is an
isomorphism, $\pi_{d-1}(\partial)$ maps the generator $\iota$ to
$\pm id_{S^d}$. Therefore up to signs, $\iota$ is a section up to
homotopy of $\partial$. By applying Lemma~\ref{splitting lemma} to
the homotopy fiber sequence $\Omega S^{2d-1}\buildrel{\Omega
f}\over\rightarrow \Omega S^d \buildrel{\partial}\over\rightarrow
S^{d-1}$, we obtain that the loop sum $\iota * \Omega
f:S^{d-1}\times\Omega S^{2d-1}\buildrel{\simeq}\over
\rightarrow\Omega S^d$ is a homotopy equivalence and in particular
$\Omega f:\Omega S^{2d-1}\rightarrow\Omega S^d$ has a retract up
to homotopy. As a consequence, by looping once again, one obtains
splittings:
\begin{eqnarray*}
\Omega^2 S^2 &\simeq& \Omega S^1 \times  \Omega^2 S^3    \\
\Omega^4 S^4 &\simeq& \Omega^3 S^3 \times \Omega^4 S^7   \\
\Omega^8 S^8 &\simeq& \Omega^7 S^7 \times \Omega^8 S^{15} \quad .
\end{eqnarray*}
And we obtain that the $4$-fold loop map $\Omega^4\nu:\Omega^4 S^7
\rightarrow  \Omega^4 S^4$ and the $8$-fold loop map
$\Omega^8\sigma:\Omega^8 S^{15}\rightarrow  \Omega^8 S^8$ are
injective in homology. In the case of $\eta$, since $\Omega^2
S^1\rightarrow\Omega^2 S^3\buildrel{\Omega^2\eta}\over \rightarrow
\Omega_0^2 S^2$ is a fibration with path connected base and
contractible fiber, we have that the two fold loop map $
\Omega^2\eta:\Omega^2 S^3\buildrel{\simeq}\over\rightarrow
\Omega_0^2 S^2 $ is a homotopy equivalence.

\m
Let
$ad_d:\pi_{i+d}(X)\buildrel{\cong}\over\rightarrow\pi_{i}(\Omega^d
X)$ denote the adjunction map. We have the following commuting
diagram
$$
\xymatrix{
\pi_{2d-1}S^{2d-1}\ar[r]^{ad_d}_{\cong}\ar[d]_{\pi_{2d-1}(f)} &
\pi_{d-1}\Omega^d S^{2d-1}\ar[r]^{\hur
}_{\cong}\ar[d]_{\pi_{d-1}(\Omega^d f)}
&H_{d-1}\Omega^d S^{2d-1}\ar[d]^{H_{d-1}(\Omega^d f)}\\
\pi_{2d-1}S^{d}\ar[r]^{ad_d}_{\cong} & \pi_{d-1}\Omega^d
S^{d}\ar[r]^{\hur }
&H_{d-1}\Omega^d S^{d}\\
& \pi_{d-1}SO(d)\ar[r]^{\hur }\ar[u]_{\pi_{d-1}(J)}\ar[lu]^{\J} &
H_{d-1}SO(d)\ar[u]_{H_{d-1}(J)} }
$$
where $\J$ is the classical $J$-homomorphism and $J$ is defined,
for example, as in \cite[section 5]{gaudensmeni0}. Recall also
that $\hur$ is the ({\it modulo } $2$) Hurewicz homomorphism.
Since $H_{d-1}(\Omega^d f)$ is injective, $\hur \circ
ad_d\circ\pi_{2d-1}(f)$ is also injective and so
$$\hur \circ ad_d\circ\pi_{2d-1}(f)(id_{S^{2d-1}})=
\hur \circ ad_d(f)\neq 0,$$ {\it i. e.} $f:=\eta$, $\nu$ or
$\sigma$ is detected by the Hurewicz homomorphism for $\Omega^d
S^d$. Since $\J(a)=\eta$ and $\hur \circ ad_2(\eta) \neq 0$, it
follows that $\hur (a)\neq 0$. Similarly $\hur (b)\neq 0$ and
$\hur (c)\neq 0$. We have finally proved:
\begin{lem}
\label{lem:hur:so} The Hopf maps $\eta\in\pi_1\Omega^2S^2$,
$\nu\in\pi_3\Omega^4S^4$, and $\sigma\in\pi_7\Omega^8S^8$ as well
as their preimages $a\in\pi_1SO(2)$, $b\in\pi_3SO(4)$,
$c\in\pi_7SO(8)$ under the $J$-homomorphism are detected by the
Hurewicz homomorphism. This holds unstably as well as stably.
\end{lem}
The stable case follows for instance  from the unstable one by
using the work of Milgram \cite{milgram}, where it is shown that
the mod $2$ homology of $SO$ injects in the mod $2$ homology
of $QS^0$ via the $J$-homomorphism.

\subsection{Homology of $\Omega^2 S^3$ and $\Omega^2 S^2$ as $BV$-algebras}

The homology of $\Omega^2 {S}^{3}$, as a Pontryagin algebra, is
polynomial on generators $u_n$ of degree $2^n -1, n\geq 1$, where
$u_n = Q_1(u_{n-1})$ (see Cohen's work in
\cite{Cohen-Lada-May:homiterloopspaces}). We first notice that
$u_1$ is the bottom non trivial class in positive degrees, and as
such, must be in the image of the Hurewicz homomorphism, according
to the Hurewicz theorem. But $\pi_1 \Omega ^2 S^3$ is infinite
cyclic generated by $\iota=\mathrm{ad}_2 (\mathrm{Id}_{S^3})$ the
adjoint of the identity of $S^3$.

\m Since $u_1=\hur \iota=\hur \circ ad_2(id_{S^3})$, thanks to
part i) of \cite[corollary 5.7]{gaudensmeni0},
$$
{BV} (u_1) = \hur  (\mathrm{ad}_2(\Sigma\eta))
$$
where $\mathrm{ad}_2(\Sigma\eta) \in \pi_2 (\Omega^2 S^3)$ is the
adjoint of suspension of the Hopf map $\Sigma\eta\in \pi_4 S^3$.
We can conclude that  $\hur  (\mathrm{ad}(\Sigma\eta))$ has to be
non zero, because the two-fold loop equivalence $\Omega^2 S^3
\longrightarrow  \Omega_0^2 S^2$ takes $\mathrm{ad} (\Sigma\eta)$
to $\Sigma\eta \circ \eta $ which is well known to be spherical
(because it is a map of Kervaire invariant one and such elements
are known to be spherical, see Remark \ref{remkerv}). We will
nevertheless give a direct argument to see that
$\mathrm{ad}(\Sigma\eta)$ is detected by the Hurewicz
homomorphism. We want to show that the Hurewicz homomorphism for
$\Omega^2 S^3$ is non trivial in degree $2$. This can be seen as
follows. We know that $\pi _1 \Omega ^2 S^3 \cong\mathbb{Z}$ and
$\pi _2 \Omega^2 S^3 \cong\mathbb{Z}/2\mathbb{Z}$ generated by
$ad_2(\Sigma \eta)$.

\sm Let $p:S^3\rightarrow  K(\mathbb{Z},3)$ represent the
generator for the third integral cohomology group of $S^3$,
$\pi_3K(\mathbb{Z},3)\cong H^3(S^3;\mathbb{Z})$. Let $S^3\langle
3\rangle$ be the homotopy fiber of $p$ and let $j:S^3\langle
3\rangle\rightarrow S^3$ be the fiber inclusion. Let
$\iota:S^1\rightarrow\Omega^2 S^3$ the adjoint of the identity of
$S^3$. Since $S^3\langle 3\rangle$ is $3$-connected,
$\pi_1(\Omega^2 p)$ is an isomorphism and so maps $\iota$ to $\pm
id_{S^1}$. Therefore by applying Lemma~\ref{splitting lemma} to
the homotopy fiber sequence $\Omega^2 S^3\langle 3\rangle
\buildrel{\Omega^2 j}\over\rightarrow \Omega^2
S^3\buildrel{\Omega^2 p} \over\rightarrow \Omega^2
K(\mathbb{Z},3)\simeq S^1$, we obtain that $\Omega^2 j:\Omega^2
S^3\langle 3\rangle\rightarrow \Omega^2 S^3$ has a retract up to
homotopy and so is injective in homology. Since
$\pi_3(S^1)=\pi_2(S^1)=0$,
 $\pi_2(\Omega^2 j):\pi_2(\Omega^2 S^3\langle 3\rangle)
\buildrel{\cong}\over\rightarrow \pi_2(\Omega^2 S^3)\cong
\mathbb{Z}/2\mathbb{Z}$ is an isomorphism. Since $\Omega^2
S^3\langle 3\rangle$ is simply connected, the Hurewicz
homomorphism $\hur :\pi_2(\Omega^2 S^3\langle 3\rangle)
\buildrel{\cong}\over\rightarrow H_2(\Omega^2 S^3\langle
3\rangle)$ is an isomorphism. Since $H_2(\Omega^2 j)$ is also an
isomorphism, the Hurewicz homomorphism
$$
\hur :\pi_2(\Omega^2
S^3)={\mathbb{Z}}/{2\mathbb{Z}}.\mathrm{ad}_2(\Sigma\eta)
\buildrel{\cong}\over\longrightarrow H_2(\Omega^2
S^3)={\mathbb{Z}}/{2\mathbb{Z}}.u_1^2
$$
is an isomorphism. Hence $\hur  (\mathrm{ad}_2(\Sigma\eta))=u_1^2$
and so $BV(u_1)=u_1^2$. We refer to the appendix  in
\cite{menibv} for a different derivation of this result.

\m As next step, we compute using Corollary \ref{comm:bv:q} that
$$
{BV}  (u_2)= {BV} (Q_1 u_1)= \{{BV} (u_1), u_1\} + {BV} (u_1)* {BV} (u_1)=  (u_1)^{4}\quad .
$$
because all Browder brackets vanish in $H_* \Omega ^2 S^3$:
$S^3$ is a Lie group, hence a loop space and therefore  $\Omega ^2
S^3$ is in fact a $3$-fold loop space (Recall that the Browder bracket is
an obstruction to extend the $2$-fold loop structure to a $3$-fold
one, see Remark \ref{rem:three:three}). We thus have ${BV}
(u_2)=(u_1)^{4}$, and more generally, for $n>1$, we get by
induction:
\begin{eqnarray*}
{BV} (u_n) &=& {BV} (Q_{1}u_{n-1}) =  \{ {BV} u_{n-1}, u_{n-1}\}+
({BV} u_{n-1})^{2} = u_1^{2^{n}}  ~~.
\end{eqnarray*}
As all brackets are trivial, ${BV}$ is a derivation for the
Pontryagin product, and we obtain:
\begin{theor}
\label{bvcalc:omegasdeux} The action of the Batalin-Vilkovisky
operator $BV$ on $H_* (\Omega^2S^3 , \mathbb{F}_2 )= \mathbb{F}_2
[u_i; i\geq 1]$ is  given on Pontryagin monomials $u_{i_1}^{\ell_1}*u_{i_2}^{\ell_2}*\ldots *u_{i_n}^{\ell_n}$ with $\ell_j>0$ by
\begin{eqnarray}
{BV} (u_1^{\ell_1}*u_2^{\ell_2}*\ldots *u_n^{\ell_n})= \sum_{j=1}^n \ell_j
u_1^{2^{i_j}}*u_{i_1}^{\ell_1}*u_{i_2}^{\ell_2}*\ldots *u_{i_j}^{\ell_j-1}* \ldots *u_{i_n}^{\ell_n}.
\end{eqnarray}
\end{theor}
where the coefficient $\ell_i$ in the sum is of course reduced mod
$2$.

Thanks to the two-fold loop equivalence $H_*\Omega_0^2 S^2\cong
H_* \Omega^2 S^3$, we can compute the $BV_2$ structure of
$H_*\Omega^2 S^2$. Indeed, let $[i]$ be the zero dimensional class
in $H_* \Omega^2 S^2$ corresponding to the degree $i$ component.
All elements of $H_* \Omega^2 X$ are uniquely of the form $x*[i]$
with $x \in H_* \Omega_0^2 S^2 $. Hence there is an isomorphism of
algebras (for the Pontryagin product)
$$
H_*\Omega^2 S^2 \cong \mathbb{F}_2 [\mathbb{Z}]\otimes H_*
\Omega^2_0 S^2.
$$
We have to compute $BV ([i])$ and Browder brackets of the form
$\{[i], u_1^{\ell_1}*u_2^{\ell_2}*\ldots *u_n^{\ell_n}\}$: the
action of $BV$ on a general element $[i]\otimes
u_1^{\ell_1}*u_2^{\ell_2}*\ldots *u_n^{\ell_n} $ follows by the
Batalin-Vilkovisky Formula (\ref{crochet defaut de BV
derivation}). We begin with the following lemma, which follows
from the discussion in \cite[section 5]{gaudensmeni0}:

\sm\begin{lem}
For any pointed topological space $X$ and for any element $g$ in
$H_* (\Omega^2 X)$, the following formula holds:
$$
BV (g) = g \circ (u_1 *[1])
$$
where $\circ$ is the map induced in homology by the composition
action:
$$
\Omega^2 S^2 \times \Omega^2 X \longrightarrow \Omega ^2 X ~~.
$$
\end{lem}
In particular, in $H_* \Omega^2 S^2$:
$$
BV ([1])= [1] \circ (u_1 * [1])=(u_1 * [1])
$$
and more generally:

\sm
\begin{lem}
For any $i\in\mathbb{Z}$, the following formula holds
$$
BV([i])= i(u_1*[i])~~.
$$
\end{lem}

Indeed,
$$
BV([2i])= BV([i]*[i])= 2 BV[i]=0
$$
while
$$
BV([2i+1])= BV([1])*[2i]= u_1 *[2i+1]
$$
by Formula (\ref{crochet defaut de BV derivation}).

\sm\begin{lem}
All Browder brackets  are trivial in $H_*\Omega^2 S^2$.
\end{lem}

{\bf Proof.~ } Firstly, all brackets are trivial on the subobject $H_* \Omega_0^2
S^2 \cong H_* \Omega^2 S^3$ as shown above. We have:
\begin{eqnarray*}
\{[1], u_1\}&=& BV ([1]*u_1)+ BV[1] * u_1+ [1]*BV (u_1) \\
& =& BV (BV [1])+ u_1*u_1*[1]+ u_1* u_1*[1]\\
&=& BV^2 ([1])+ 2  (u_1*u_1*[1]) \\
&=& 0 ~~.
\end{eqnarray*}
By induction on $i$, the brackets $\{[i], u_1\}$ are trivial for
all $i\in \mathbb{Z}$.  Now, there is a general formula in the
homology of two fold loop spaces \cite[Theorem 1.3 (4) p.
218]{Cohen-Lada-May:homiterloopspaces}
\begin{eqnarray}
\label{formulll} \{x, Q_1 y\}= \{\{x, y\}, y\}
\end{eqnarray}
hence
$$
\{[1], u_2\}= \{[1], Q_1 u_1\}= \{\{[1], u_1\}, u_1\}=0
$$
and more generally, using Formula (\ref{formulll}) inductively we
get for $u_{n+1}=Q_1 u_n$ that $\{[1], u_{n+1}\}=0$. We conclude
with the Formula  (\ref{crochet defaut de BV derivation}) that all
brackets are trivial in $H^* \Omega^2 S^2$.
\hfill $\square$

Gathering all these calculations, we obtain:

\sm\begin{theor}
\label{omegadeuxsdeux} The $BV$ structure on $H_* \Omega_2 S^2
\cong \mathbb{F}_2 [\mathbb{Z}]\otimes H_* \Omega^2 S^3$ is given
by restriction on $H_* \Omega^2 S^3$ (theorem above) by the on the
component $[0]\otimes H_* \Omega^2 S^3$. For $i\neq 0$ and $f\in
H_* \Omega ^2 S^3$, the action of the $BV$ operator is given by:
$$
BV ([i]* f) = BV([i])*f + [i]* BV(f)~
$$
with  $BV([i])=[i]*u_1$ if $i$ is odd and zero otherwise.
\end{theor}

\subsection{The $BV$-structure of $H_* \Omega^2 S^{k+2}$ for $k>1$}

The homology of $\Omega^2 S^{k +2}$ is, as a Pontryagin algebra,
the polynomial algebra $\mathbb{F}_2[u_n,n\geq 1]$ on generators
$u_n=Q_1u_{n-1}$ of degree $2^{n-1}(k+1)-1$ (see Cohen's work in
\cite{Cohen-Lada-May:homiterloopspaces}). By
Corollary~\ref{comm:bv:q}, for  $n>1$,
\begin{eqnarray*}
{BV} (u_n) = {BV} (Q_1 u_{n-1})=  \{{BV} (u_{n-1}), u_{n-1}\}.
\end{eqnarray*}
Since $k>1$, ${BV} (u_1) =0$, for degree reasons. Hence by
induction, ${BV} (u_n)=0$ for all $n\geq 1$.

We are going now to see that all the Browder brackets vanish in
$H_*\Omega^2 S^{k+2}$. According to~\cite[Theorem 1.2 (3) p.
215]{Cohen-Lada-May:homiterloopspaces}, $\{x, x\}=0$, and
according to the Formula (\ref{formulll}). Therefore
$\{u_n,u_n\}=0$ and
$\{u_{n+i+1},u_n\}=\{u_{n+i},\{u_{n+i},u_n\}\}$. So by induction,
$\{u_{n+i},u_n\}=0$. By anti-commutativity of the Browder bracket,
for all $m,n\geq 1$, $\{u_m,u_n\}=0$. Since the Browder brackets
are trivial on the generators, using the Poisson relation, all the
Browder brackets are trivial on $H_* \Omega^2 S^{k+2}$. Therefore
by~(\ref{crochet defaut de BV derivation}), the operator $BV$ is a
derivation with respect to the Pontryagin product. It follows that
$BV$ is the trivial operator.

\sm\begin{theor}
\label{trivialahah} The operator $BV$ acts trivially on the {\it
modulo} $2$ homology of $\Omega^{2}S^{2+k}$ for $k>1$.
\end{theor}


\section{Some calculations of the {\it mod} $2$ Hurewicz homomorphism for $QS^0$}
\label{comp:hur:qso}

We wish to describe a method to get information on the {\it mod}
$2$ Hurewicz homomorphism for $QS^0$. We recall that this map is
conjecturally described by the conjecture stated
\cite{minamikervaire} under the name of \emph{Curtis-Madsen}
conjecture.
\begin{conjecture}~\cite[chap. 1]{minamikervaire}
\label{curtismadsen}
Let $f\in\pi_{>0} QS^0$ be any stable map of positive degree. Then
$f$ is not in the kernel of the {\it mod} $2$ Hurewicz
homomorphism if and only if $f$ has Hopf or Kervaire invariant
one.
\end{conjecture}
\begin{rem}
\label{remkerv} It is well known, in connection with the Kervaire
invariant problem, that maps of Kervaire invariant one (when they
do exist) are detected by the Hurewicz homomorphism (see for instance \cite{madsenmilg}, \cite[p. 30]{snaith}). The
composition squares of the Hopf maps are maps of Kervaire
invariant one, and we will recover in the following that they are
detected by the Hurewicz homomorphism by direct calculations.
\end{rem}
There is for each $n>0$ a natural inclusion $SO(n)\to \Omega^n
S^n$ by compactifying the natural action of $SO(n)$ on
$\mathbb{R}^n$. These inclusions are compatible as $n$ grows, and
stabilize to a map  $\Theta: SO\to Q_1 S^0$. We notice that the
effect of
\begin{itemize}
\item the $\circ$-product on $H_* Q_1S^0$, \item the map $H_*
\Theta: H_* SO \longrightarrow   H_* Q_1S^0$, \item  and the map
$J: \pi _*SO \longrightarrow   \pi_* QS^0$,
\end{itemize}
are completely known by \cite{madsen,milgram}. Let us only recall
what is relevant for our purposes. First the map
$H_*\Theta:H_*SO\rightarrow H_*Q_1S^0$ is injective and takes the
Pontryagin product in $H_*SO$ to the circle product (the product
in homology  induced by the topological monoid structure given by
composition of loops). Finally we point out that $H_* QS^0$ is a
polynomial algebra for the Pontryagin loop sum.

\m We notice that we already know that the following elements
$\eta$, $\nu$, $\sigma$, and $\eta ^2$ are spherical detected by
the Hurewicz homomorphism. Once again, it possible to provide a
direct argument for this fact in the special case of $\nu^2$ and
$\sigma^2$. Let us check the case of $\nu ^2$, as the case of
$\sigma^2$ is quite similar. One has:
\begin{eqnarray}
\hur  (\nu^2)&=& \hur  (\nu). \hur (\nu)\label{uno}\\
&=& \Delta_{p_3} (\hur (\nu))\label{dos}\\
&=& \Delta_{p_3} (\theta *[-1])\label{tres}
\end{eqnarray}
where $\theta= H_* (\Theta) (\hur  (b))$. The equality (\ref{uno})
expresses the fact that the action of $SO(n)$ on loops spaces
corresponds to the composition (see  \cite[Theorem 5.6]{gaudensmeni0}). The equality
(\ref{dos}) corresponds to the fact that $\hur (\nu)$ is non zero
and must be primitive  in $H_* SO(n)$ (as anything lying in the
image of $\mathcal{H}$), hence its action is that of
$\Delta_{p_3}$. Finally the equation (\ref{tres}) is the
observation that $\nu= J(a)$ is in the image of the
$J$-homomorphism, that decomposes as $\J= (-*[-1])\Theta$.
$\Delta_{p_3}$ is a derivation for the Pontryagin product, and we
have
$$
\Delta_{p_3} (\theta *[-1])= \Delta_{p_3} (\theta)*[-1] + \theta
* \Delta_{p_3} [-1] \quad .
$$
But $\Delta _{p_3}$ is the composition product with $\theta$, hence
$$
\Delta_{p_3} (\theta *[-1])= (\theta \circ \theta)*[-1]+ \theta
*(\theta \circ [-1]) = \theta* \chi \theta = \theta *\theta *[-2]=
(\theta*[-1])^{2}
$$
because $H_* SO$ is an exterior algebra. Here $\chi$ is the
antipode of the Hopf algebra $H_* QS^0$. In particular, $\hur (\nu^2)$ is non zero. One shows in the same way that
$\sigma^2
$ is detected by the Hurewicz
homomorphism, consistently with Remark \ref{remkerv}.

Another observation that one can do is:
\begin{eqnarray*}
\hur  (\eta^3)= \hur  (\eta). \hur (\eta ^2)= \Delta_1 (\hur
(\eta).\hur (\eta) )= (BV)^2 (\hur  (\eta))=0 \quad .
\end{eqnarray*}
Because ${BV}$ has order $2$, as well as any higher $BV$-operator because $H_* SO(n)$ is an exterior algebra. This work in the same way for $\nu$
and $\sigma$, and we deduce:
\begin{eqnarray*}
\hur (\nu^3) = \hur  (\sigma^3)=0  \quad .
\end{eqnarray*}
Alternatively, one could argue as follows to show for instance that $\hur (\eta^3)$ vanishes. The element $\hur
(\eta)$ is $J_* (\hur  (a))$ where $a$ was the generator of $\pi_1
SO$. Hence $\hur  (a)$ is a homogeneous polynomial of odd degree
in $H_* SO$. But from Corollary \ref{cor:nul:sur:square}, we know
that the action of odd degree homogeneous elements in $H_* SO$ on
Pontryagin squares is trivial, hence:
$$
\hur (\eta^3)= \Delta_1 (\hur (\eta^2))= \Delta_1
(\theta*[-1])^{2})=0 ~.
$$
This works all the same for $\nu^3$ and $\sigma^3$.

To go further, we need one more result, which the reader might
draw from the  information in \cite{semc}.

\sm\begin{proposition}
\label{hurewicz:so} The {\it mod} $2$ Hurewicz homomorphism for
$\mathrm{SO}$ annihilates all classes except those yielding an
element of Hopf invariant one.
\end{proposition}
We note that we have already proved (Lemma \ref{lem:hur:so}) that
the classes $a$, $b$, and $c$ are detected by the Hurewicz
homomorphism for $SO$. The fact that the rest of the homotopy is
not follows from a simple computation using Bott periodicity and
the information in \cite{semc}.

\sm By a result of Novikov \cite{novikov} if $\alpha $ is in $(\mathrm{Im}
J)_i$ and $\beta$ is in $\Theta_j$ with $j<2i$ then $\alpha \beta$
is again in $\mathrm{Im} J$. Here $\Theta_j$ classically denotes
the subgroup of the $i^{\mathrm{th}}$ stable homotopy group $\pi_i
^S$ of the sphere spectrum that consists in elements representable
by homotopy spheres. The classical works on exotic spheres show
that the quotient of the stable stem $\pi_i^S$ by $\Theta_i$ has
at most order two, and non triviality of the cokernel happens only
in degrees of the form $2^n-2$.  Hence from the above Proposition
\ref{hurewicz:so}, we could already deduce particular cases of
Theorem \ref{divtrivial}, which is however much more general.
\begin{defin} We say that an element $\theta \in \pi_* QS^0 \cong
\pi_*^S$ is divisible by $x \in \pi_*QS^0 \cong \pi_*^S$ if
$\theta=x
\theta'$ for some $\theta'$.
\end{defin}
Hence the set of elements divisible by those of a fixed subset of $\pi_*^S$ is the ideal generated by this subset.
Now we can state a result that fits perfectly with the
Curtis-Madsen conjecture~recalled in~\ref{curtismadsen}:

\sm\begin{theor}
\label{divtrivial} All classes in $\pi_{>0} QS^0$ divisible by the
image of the $J$-homomorphism are annihilated by  the {\it mod}
$2$ Hurewicz homomorphism, except the Hopf maps and their
composition squares which are not annihilated. In other words, Conjecture \ref{curtismadsen} holds on the ideal generated by the image of the $J$-homomorphism.
\end{theor}
Before proceeding to the proof, one might wonder how strong this
result is. Are there many elements that are not in the image of
the $J$-homomorphism, but divisible by the image of the
$J$-homomorphism? The general problem of describing the behaviour
of the multiplication by elements in the image of the
$J$-homomorphism in $\pi_*^S$ is probably very difficult.
Nevertheless,  one can explicitly construct some infinite families
of elements in the stable homotopy groups of spheres that are
divisible by the image of the $J$ homomorphism, but are not in the
image of the $J$-homomorphism. This can be done for instance by
using \emph{the spectrum of topological modular forms}. In
\cite{bauercomput}, a computation of the homotopy groups of this
spectrum is produced. We see that
\begin{itemize}
\item the spectrum of topological modular forms $\mathit{tmf}$ is
a ring spectrum, \item the unit map $S^0\longrightarrow
\mathit{tmf}$ detects the Hopf maps $\eta$ and $\nu$, \item the
homotopy groups of $\mathit{tmf}$ have a periodicity of order
$192$, \item there are periodic families of chromatic filtration
$2$ that are detected by $S^0\longrightarrow \mathit{tmf}$, and
whose image support a non trivial multiplication by $\eta$ or
$\nu$ in the homotopy of $\mathit{tmf}$.
\end{itemize}
So, each of these families provides infinitely many elements
divisible by the image of the $J$-homomorphism but not in the
image of the $J$-homomorphism. A concrete example is given by the
periodic family generated by $\eta\kappa$, which is divisible by
$\eta$ (see \cite[p. 24]{mahohopell}).

\sm
{\bf Proof.~}As we noticed before, the Hopf maps $\eta$, $\nu$ and $\sigma$,
and their composition squares are detected by the Hurewicz
homomorphism (Lemma \ref{lem:hur:so}), while for $n>2$, $\eta^n$,
$\nu^n$ and $\sigma^n$ are all annihilated by the Hurewicz
homomorphism (observe that in fact these elements are zero for $n>0$). Assume now that $\varphi=\psi
\theta$,
where $\psi$ is an element of the image of the $J$-homomorphism
which is not a Hopf map. Then we have
\begin{eqnarray*}
\hur (\varphi)= \hur (\psi
\theta)=  \hur (\psi)\circ
\hur  (\theta)~.
\end{eqnarray*}
But Proposition \ref{hurewicz:so} states that $\hur (\psi)$ is
trivial if  $\psi$ is not a Hopf map. It remains to prove that  $
\hur (\psi
\theta)=0$ as soon as $\psi$ is one of the
Hopf maps. The first case is $\psi= \eta$. We have:
$$
\hur  (\eta
\theta)= \hur   (\eta ). \hur  (\theta)=
\Delta _1 (\hur (\theta)) \quad .
$$
This follows from \cite{gaudensmeni0}
Recall that $H_* QS^0$ is generated by Kudo-Araki operations
applied to the element $[1]$ and shifting components by $-*[i]$
for all $i$. In particular $\hur  (\eta)$ being non zero in $H_1
Q_0S^0 $, we have:
$$
\hur  (\eta) = Q_1 ([1])* [-2]= Q^1([1])* [-2]
$$
Indeed, the only non zero element of degree one in $H_* QS^0$ are
of the form $Q^1([1])* [i]$, $Q^1([1])$ belongs to $H_* Q_2 S^0 =
H_* Q_0S^0 *[2] \subset H^* QS^0$, and the Hurewicz homomorphism
lands in $H^* Q_0 S^0\subset H_* QS^0$ as $Q_0 S^0$ is the
component of the basepoint.

\sm According to \cite[theorem 6.18]{madsenmilg}, there is a
formula
\begin{eqnarray}
\label{formule} Q^i [1] \circ a = \sum _k Q^{i+k} (Sq_*^k a)
\end{eqnarray}
where  $Sq_*^k$ are dual to the Steenrod operations. In our case,
we obtain:
\begin{eqnarray*}
\hur  (\eta
\theta)= \sum _k Q^{1+k} Sq_*^k (\hur
(\theta))~.
\end{eqnarray*}
On the other hand, all Steenrod operations of strictly positive
degree vanish on spherical elements, hence $\hur
(\eta
\theta)=  Q^{1} (\hur (\theta))$. We conclude that
by instability (see Remark \ref{renorm}) that $\hur
(\eta
\theta)=  0$ if $\theta$ has degree strictly greater
than $1$. On the other hand in degree $1$, the only possible non
trivial element is $\theta=\eta$, and we already know that $\eta
\eta$ is detected by the Hurewicz homomorphism.

This method applies as well for $\nu$ and $\sigma$. Wellington
\cite[remark 5.8, p. 47]{wellington} provides formula for $\hur
(\nu)$ and $\hur  (\sigma)$ in terms of Kudo-Araki operations
applied to $[1]$.  Formula (\ref{formule}) generalizes to the
formula (proposition 1.6 of May's article \emph{The homology of
$E_\infty$-ring spaces} in
\cite{Cohen-Lada-May:homiterloopspaces}):
\begin{eqnarray}
\label{actionhj} Q^i x \circ a = \sum _k Q^{i+k}(x Sq_*^k a)\quad
.
\end{eqnarray}
This, together with an instability argument yields:

\sm\begin{lem}
If the degree of $\theta$ is strictly bigger than $7$, $\hur
(\nu
\theta)$ and $\hur (\sigma
\theta)$ are
trivial.
\end{lem}
In degree less than $7$, a system of generators for $\pi_*^S$ is
given by $\{\eta, \eta^2, \nu, \nu^2, \sigma\}$. These cases have
already been taken care of. This finishes the proof.

\appendix

\section{Two definitions}
\label{previous}

We recall first the definition of Gerstenhaber algebras.
\sm\begin{defin}
\label{definition algebre de Gerstenhaber} A {\it $e_n$-algebra}
is a commutative graded algebra $A$ equipped with a linear map
$\{-,-\}:A \otimes A \to A$ of degree $n-1$ such that:

\noindent a) the bracket $\{-,-\}$ gives $A$ the structure of
graded Lie algebra of degree $n-1$. This means that for each $a$,
$b$ and $c\in A$

$\{a,b\}=-(-1)^{(\vert a\vert+n-1)(\vert b\vert+n-1)}\{b,a\}$ and

$\{a,\{b,c\}\}=\{\{a,b\},c\}+(-1)^{(\vert a\vert+n-1)(\vert
b\vert+n-1)} \{b,\{a,c\}\}.$

\noindent b)  the product and the Lie bracket satisfy the Poisson
relation:
$$\{a,bc\}=\{a,b\}c+(-1)^{(\vert a\vert+n-1)\vert b\vert}b\{a,c\}.$$
\end{defin}
We also recall the definition of a  $BV_n$-algebra.
\sm\begin{defin}~\cite[Def 5.2]{Salvatore-Wahl:FrameddoBVa}
\label{definition BV_n-algebre}
A $BV_n$-algebra $A$ is an $e_n$-algebra with a linear
endomorphism $BV:A\rightarrow A$ of degree $n-1$ such that
$BV\circ BV=0$ and for each $a,b\in A$,
\begin{equation}\label{crochet defaut de BV derivation}
\{a,b\}=(-1)^{\vert a\vert}\left(BV(ab)-(BVa)b-(-1)^{\vert a\vert}
a(BVb)\right).
\end{equation}
The bracket measures the deviation of the operator $BV$ from being
a derivation with respect to the product.
\end{defin}


\bibliography{Bibliographie_Gerald}
\bibliographystyle{amsplain}

\end{document}